\documentclass[a4paper, reqno]{amsart}

\title[{Finite $\mathcal{W}$-algebras of $\mathfrak{osp}_{1|2n}$ and Ghost centers}]{Finite $\mathcal{W}$-algebras of $\mathfrak{osp}_{1|2n}$ and Ghost centers}

\author{Naoki Genra}
\address{Kavli IPMU (WPI), UTIAS, The University of Tokyo, Kashiwa, Chiba 277-8583, Japan}
\email{naoki.genra@ipmu.jp}

\usepackage{latexsym}
\usepackage{amsmath}
\usepackage{amsfonts}
\usepackage{amssymb}
\usepackage{amsxtra}
\usepackage{mathrsfs}
\usepackage[all]{xy}

\newtheorem{definition}{Definition}[section]
\newtheorem{proposition}[definition]{Proposition}
\newtheorem{theorem}[definition]{Theorem}
\newtheorem{corollary}[definition]{Corollary}
\newtheorem{lemma}[definition]{Lemma}
\newtheorem{remark}[definition]{Remark}

\newtheorem*{thmA}{Theorem A}

\numberwithin{equation}{section}

\newcommand{\Z}{\mathbb{Z}}

\newcommand{\R}{\mathbb{R}}
\newcommand{\C}{\mathbb{C}}
\newcommand{\End}{\operatorname{End}}

\newcommand{\Span}{\operatorname{Span}}
\newcommand{\id}{\operatorname{id}}

\newcommand{\NO}[1]{:\!#1\!:}
\newcommand{\vac}{|0\rangle}

\begin{document}
\maketitle
\markboth{Finite $\mathcal{W}$-algebras of $\mathfrak{osp}_{1|2n}$ and Ghost Centers}{Finite $\mathcal{W}$-algebras of $\mathfrak{osp}_{1|2n}$ and Ghost Centers}
\begin{abstract}
We prove that the finite $\mathcal{W}$-algebra associated to $\mathfrak{osp}_{1|2n}$ and its principal nilpotent element is isomorphic to Gorelik's ghost center of $\mathfrak{osp}_{1|2n}$, which proves an analog of Kostant's theorem for $\mathfrak{osp}_{1|2n}$.
\end{abstract}

\section{Introduction}

A Lie suparalgebra $\mathfrak{osp}_{1|2n}$ is a finite-dimensional simple Lie superalgebra whose Dynkin diagram is the same as type $B_n$ except for a unique simple short root, which is replaced by a non-isotropic odd simple root in $\mathfrak{osp}_{1|2n}$. The Lie suparalgebra $\mathfrak{osp}_{1|2n}$ is not a Lie algebra but has similar properties to simple Lie algebras. For example, the category of finite-dimensional $\mathfrak{osp}_{1|2n}$-modules is semisimple and we have the Harish-Chandra isomorphism $Z(\mathfrak{osp}_{1|2n}) \simeq \C[\mathfrak{h}]^W$, where $Z(\mathfrak{g})$ denotes the center of the universal enveloping algebra $U(\mathfrak{g})$, $\mathfrak{h}$ is a Cartan subalgebra of $\mathfrak{osp}_{1|2n}$ and $W$ is the Weyl group. However $\mathfrak{osp}_{1|2n}$ doesn't satisfy the Duflo theorem \cite{Duflo77}, which says that annihilators of Verma modules in $U(\mathfrak{g})$ is generated by its intersections with the center $Z(\mathfrak{g})$ for simple Lie algebras $\mathfrak{g}$. This problem was founded by Musson \cite{Musson97} and solved by Gorelik and Lantzmann \cite{GL} by using an extension algebra $\widetilde{Z}(\mathfrak{osp}_{1|2n})$ of $Z(\mathfrak{osp}_{1|2n})$. More precisely, Gorelik and Lantzmann prove that annihilators of Verma modules in $U(\mathfrak{osp}_{1|2n})$ is generated by its intersections with $\widetilde{Z}(\mathfrak{osp}_{1|2n})$. The associative algebra $\widetilde{Z}(\mathfrak{osp}_{1|2n})$ is called the ghost center of $\mathfrak{osp}_{1|2n}$ in \cite{Gorelik01}.

For a Lie superalgebra $\mathfrak{g}$ with $\mathfrak{g}_{\bar{1}} \neq0$, the ghost center $\widetilde{Z}(\mathfrak{g})$ is introduced by Gorelik in \cite{Gorelik01} as the direct sum $Z(\mathfrak{g}) \oplus \mathcal{A}(\mathfrak{g})$, where $\mathcal{A}(\mathfrak{g})$ is the anticenter defined by $
\mathcal{A}(\mathfrak{g}) = \{a \in U(\mathfrak{g}) \mid ua-(-1)^{p(u)(p(a)+\bar{1})}au = 0\ \mathrm{for}\ \mathrm{all}\ u \in \mathfrak{g}\}$. If $\mathfrak{g}$ is a finite-dimensional simple basic classical Lie superalgebra, it is known that $\widetilde{Z}(\mathfrak{g})$ coincides with the center of $U(\mathfrak{g})_{\bar{0}}$ and thus is a purely even subalgebra of $U(\mathfrak{g})$. Moreover, if $\mathfrak{g}=\mathfrak{osp}_{1|2n}$, there exists $T \in U(\mathfrak{g})_{\bar{0}}$ such that $\mathcal{A}(\mathfrak{osp}_{1|2n}) = Z(\mathfrak{osp}_{1|2n}) T$ by \cite{ABF, Musson97, GL}. The element $T$ is called the Casimir ghost \cite{ABF} since $T^2 \in Z(\mathfrak{osp}_{1|2n})$. In case $\mathfrak{g}=\mathfrak{osp}_{1|2}$, \cite{Pinczon} also suggested that $T = 4Q - 4C +\frac{1}{2}$ satisfies $T^2 = 4C + \frac{1}{4}\in Z(\mathfrak{osp}_{1|2})$, where $C$ is the Casimir element in $U(\mathfrak{osp}_{1|2})$ and $Q$ is one in $U(\mathfrak{sl}_2)$.

The finite $\mathcal{W}$-algebra $U(\mathfrak{g}, f)$ is an associative superalgebra over $\C$ defined from a simple basic classical Lie superalgebra $\mathfrak{g}$ and its even nilpotent element $f$ \cite{Premet02, Premet07, Kostant, Lynch, BT, RaSo, GG}. If $\mathfrak{g}$ is a simple Lie algebra and $f$ is a principal nilpotent element $f_\mathfrak{prin}$, the corresponding finite $\mathcal{W}$-algebra $U(\mathfrak{g}, f_\mathrm{prin})$ is isomorphic to the center $Z(\mathfrak{g})$ of $U(\mathfrak{g})$ by Kostant \cite{Kostant}.

The $\mathcal{W}$-algebra $\mathcal{W}^k(\mathfrak{g}, f)$ is a vertex superalgebra defined by the Drinfeld-Sokolov reductions associated to $\mathfrak{g}, f$ and level $k \in \C$ \cite{FF92, KRW}. In general, (Ramond-twisted) simple modules of a $\frac{1}{2}\Z$-graded vertex superalgebras $V$ with a Hamiltonian operator $H$ are classified by the associated superalgebra named as the ($H$-twisted) Zhu algebras of $V$. De Sole and Kac shows that the $H$-twisted Zhu algebra of $\mathcal{W}^k(\mathfrak{g}, f)$ is isomorphic to the finite $\mathcal{W}$-algebra $U(\mathfrak{g}, f)$. In particular, there exists a one-to-one correspondence between simple modules of $U(\mathfrak{g}, f)$ and Ramond-twisted simple positive-energy modules of $\mathcal{W}^k(\mathfrak{g}, f)$. If $f = f_\mathrm{prin}$, the corresponding $\mathcal{W}$-algebra is called the principal $\mathcal{W}$-algebra of $\mathfrak{g}$, which we denoted by $\mathcal{W}^k(\mathfrak{g}) = \mathcal{W}^k(\mathfrak{g}, f_\mathrm{prin})$.

\begin{thmA}[Theorem \ref{thm:ghostZ_finiteW}]
$U(\mathfrak{osp}_{1|2n}, f_\mathrm{prin})$ is isomorphic to $\widetilde{Z}(\mathfrak{osp}_{1|2n})$ as associative algebras.
\end{thmA}

The finite $\mathcal{W}$-algebra $U(\mathfrak{osp}_{1|2n}, f_\mathrm{prin})$ associated to $\mathfrak{osp}_{1|2n}$ and its principal nilpotent element $f_\mathrm{prin}$ is an associative superalgebra with its non-trivial odd part, while the ghost center $\widetilde{Z}(\mathfrak{osp}_{1|2n})$ is not. However, we prove an isomorphism between them. Through the isomorphism in Theorem A, a $\Z_2$-grading of $\widetilde{Z}(\mathfrak{osp}_{1|2n})$ is inherited from one of $U(\mathfrak{osp}_{1|2n}, f_\mathrm{prin})$ so that the even part of $\widetilde{Z}(\mathfrak{osp}_{1|2n})$ is $Z(\mathfrak{osp}_{1|2n})$ and the odd part is $\mathcal{A}(\mathfrak{osp}_{1|2n})$.

To prove Theorem A, we use the Miura map $\mu$ and its injectivity and relationship with the Harish-Chandra homomorphism of $\mathfrak{osp}_{1|2n}$. See Section \ref{sec:finite-W} for the definition of $\mu$. The map $\mu$ was originally introduced in \cite{Lynch}. The injectivity of $\mu$ was only known for non-super cases, but has been recently proved by \cite{Nakatsuka} for super cases. As a corollary of Theorem A, it follows that simple positive-energy Ramond-twisted modules of principal $\mathcal{W}$-algebras $\mathcal{W}^k(\mathfrak{osp}_{1|2n})$ are classified by simple modules of the ghost center of $\mathfrak{osp}_{1|2n}$. See also Corollary \ref{cor:Zhu_twisted}. We remark that our definitions of $U(\mathfrak{osp}_{1|2n}, f_\mathrm{prin})$ is different from those in some literatures \cite{Poletaeva13, PS, ZS}. See Remark \ref{rem:def-W}.
\smallskip

The paper is organized as follows. In Sect.\ref{sec:Gamma/Z}, we introduce $H$-twisted Zhu algebras. In Sect.\ref{sec:W-alg}, we recall the definitions of $\mathcal{W}$-algebras $\mathcal{W}^k(\mathfrak{g}, f)$. In Sect.\ref{sec:finite-W}, we introduce two definitions $U(\mathfrak{g}, f)_{I}$ and $U(\mathfrak{g}, f)_{I\hspace{-.2em}I}$ of finite $\mathcal{W}$-algebras and show the equivalence of the definitions, that is, $U(\mathfrak{g}, f)_{I} \simeq U(\mathfrak{g}, f)_{I\hspace{-.2em}I}$. The proof is similar to \cite{D3HK}. In Sect.\ref{sec:prin-W}, we recall the principal $\mathcal{W}$-algebras $\mathcal{W}^k(\mathfrak{osp}_{1|2n})$ of $\mathfrak{osp}_{1|2n}$. In Sect.\ref{sec:Zhu_prinW}, we prove Theorem A.

\vspace{3mm}

{\it Acknowledgments}\quad The author wishes to thank Thomas Creutzig, Tomoyuki Arakawa, Hiroshi Yamauchi and Maria Gorelik for valuable comments and suggestions.  Some part of this work was done while he was visiting 
Instituto de Matem\'{a}tica Pura e Aplicada, Brazil in March and April 2022 and the Centre de Recherches Math\'{m}atiques, Universit\'{e} de Montr\'{e}al, Canada in October 2022. He is grateful to those institutes for their hospitality. He is supported by World Premier International Research Center Initiative (WPI), MEXT, Japan and JSPS KAKENHI Grant Number JP21K20317.

\section{$H$-twisted Zhu algebras}\label{sec:Gamma/Z}

Let $V$ be a vertex superalgebra. Denote by $\vac$ the vacuum vector, by $\partial$ the translation operator, by $p(A)$ the parity of $A \in V$ and by $Y(A, z) = A(z) = \sum_{n \in \Z}A_{(n)}z^{-n-1}$ the field on $V$ corresponding to $A \in V$. Let
\begin{align*}
[A_\lambda B] = \sum_{n=0}^\infty \frac{\lambda^n}{n!}A_{(n)}B \in \C[\lambda]\otimes V
\end{align*}
be the $\lambda$-bracket of $A$ and $B$ for $A, B \in V$. A Hamiltonian operator $H$ on $V$ is a semisimple operator on $V$ satisfying that $[H, Y(A, z)] = z\partial_zY(A, z) + Y(H(A), z)$ for all $A \in V$. The eigenvalue of $H$ is called the conformal weight. If $V$ is conformal and $L(z) = \sum_{n \in \Z}L_nz^{-n-2}$ is the field corresponding to the conformal vector of $V$, we may choose $H = L_0$ as the Hamiltonian operator.

Suppose that $V$ is a $\frac{1}{2}\Z$-graded vertex superalgebra with respect to a Hamiltonian operator $H$. Denote by $\Delta_A$ the conformal weight of $A \in V$. Define the $*$-product and $\circ$-product of $V$ by
\begin{align*}
A * B =\sum_{j=0}^\infty\binom{\Delta_A}{j}A_{(j-1)}B,\quad
A \circ B =\sum_{j=0}^\infty\binom{\Delta_A}{j}A_{(j-2)}B,\quad
A, B \in V.
\end{align*}
Then the quotient space
\begin{align*}
\operatorname{Zhu}_H V = V/V\circ V
\end{align*}
has a structure of associative superalgebra with respect to the product induced from $*$, and is called the $H$-twisted Zhu algebra of $V$. Here $V \circ V = \Span_\C\{A \circ B \mid A, B \in V\}$. The vacuum vector $|0\rangle$ defines a unit of $\operatorname{Zhu}_H V$. A superspace $M$ is called a Ramond-twisted $V$-module if $M$ is equipped with a parity-preserving linear map
\begin{align*}
Y_M \colon M \ni A \rightarrow Y_M(A, z) = \sum_{n \in \Z + \Delta_A}A^M_{(n)}z^{-n-1} \in (\End M)[\![z^{\frac{1}{2}}, z^{-\frac{1}{2}}]\!]
\end{align*}
such that (1) for each $C \in M$, $A^M_{(n)}C = 0$ if $n \gg 0$, (2) $Y_M(|0\rangle, z) = \id_M$ and (3) for any $A, B \in V$, $C \in M$, $n \in \Z$, $m \in \Z + \Delta_A$ and $\ell \in \Z + \Delta_B$, 
\begin{align*}
&\sum_{j=0}^\infty(-1)^j\binom{n}{j}\left(
A^M_{(m+n-j)}B^M_{(\ell+j)} - (-1)^{p(A)p(B)}B^M_{(\ell+n-j)}A^M_{(m+j)}
\right)C\\
&= \sum_{j=0}^\infty \binom{m}{j}\left(A_{(n+j)}B
\right)^M_{(m+\ell-j)}C.
\end{align*}
Hence the Ramond-twisted module is a twisted module of $V$ for the automorphism $\mathrm{e}^{2\pi i H}$. In particular, $M$ is just a $V$-module if $V$ is $\Z$-graded. Define $A^M_n$ by $Y_M(A, z) = \sum_{n \in \Z}A^M_n z^{-n-\Delta_A}$ for $A \in V$. A Ramond-twisted $V$-module $M$ is called positive-energy if $M$ has an $\R$-grading $M = \bigoplus_{j \in \R}M_j$ with $M_0 \neq 0$ such that $A^M_n M_j \subset M_{j+n}$ for all $A \in V$, $n \in \Z$ and $j \in \R$. Then $M_0$ is called the top space. By \cite[Lemma 2.22]{DK}, a linear map $V_\Gamma \ni A \mapsto A^M_0|_{M_0} \in \End M_0$ induces a homomorphism $\operatorname{Zhu}_H V \rightarrow \End M_0$. Thus we have a functor $M \mapsto M_0$ from the category of positive-energy Ramond-twisted $V$-modules to the category of $\Z_2$-graded $\operatorname{Zhu}_H V$-modules. By \cite[Theorem 2.30]{DK}, these functors establish a bijection (up to isomorphisms) between simple positive-energy Ramond-twisted $V$-modules and simple $\Z_2$-graded $\operatorname{Zhu}_H V$-modules.

\section{$\mathcal{W}$-algebras}\label{sec:W-alg}

Let $\mathfrak{g}$ be a finite-dimensional simple Lie superalgebra with the normalized even supersymmetric invariant bilinear form $(\cdot|\cdot)$ and $f$ be a nilpotent element in the even part of $\mathfrak{g}$. Then there exists a $\frac{1}{2}\Z$-grading on $\mathfrak{g}$ that is good for $f$. See \cite{KRW} for the definitions of good gradings and \cite{EK, Hoyt} for the classifications. Let $\mathfrak{g}_j$ be the homogeneous subspace of $\mathfrak{g}$ with degree $j$. The good grading $\mathfrak{g}=\bigoplus_{j\in\frac{1}{2}\Z}\mathfrak{g}_j$ for $f$ on $\mathfrak{g}$ satisfies the following properties:
\begin{enumerate}
\item $[\mathfrak{g}_i, \mathfrak{g}_j]\subset\mathfrak{g}_{i+j}$,
\item $f\in\mathfrak{g}_{-1}$,
\item $\operatorname{ad}(f)\colon\mathfrak{g}_j\rightarrow\mathfrak{g}_{j-1}$ is injective for $j\geq\frac{1}{2}$ and surjective for $j\leq\frac{1}{2}$,
\item $(\mathfrak{g}_i|\mathfrak{g}_j)=0$ if $i+j\neq 0$,
\item $\dim\mathfrak{g}^f=\dim\mathfrak{g}_0+\dim\mathfrak{g}_{\frac{1}{2}}$, where $\mathfrak{g}^f$ is the centralizer of $f$ in $\mathfrak{g}$.
\end{enumerate}
Then we can choose a set of simple roots $\Pi$ of $\mathfrak{g}$ for a Cartan subalgebra $\mathfrak{h}\subset\mathfrak{g}_0$ such that all positive root vectors lie in $\mathfrak{g}_{\geq0}$. Denote by $\Delta_j = \{ \alpha \in \Delta \mid \mathfrak{g}_\alpha \subset \mathfrak{g}_j\}$ and $\Pi_j = \Pi\cap\Delta_j$ for $j \in \frac{1}{2}\Z$. We have $\Pi=\Pi_0\sqcup\Pi_{\frac{1}{2}}\sqcup\Pi_1$. Let $\chi\colon\mathfrak{g}\rightarrow\C$ be a linear map defined by $\chi(u)=(f|u)$. Since $\operatorname{ad}(f)\colon\mathfrak{g}_{\frac{1}{2}}\rightarrow\mathfrak{g}_{-\frac{1}{2}}$ is an isomorphism of vector spaces, the super skew-symmetric bilinear form $\mathfrak{g}_{\frac{1}{2}}\times\mathfrak{g}_{\frac{1}{2}}\ni(u, v)\mapsto\chi([u, v])\in\C$ is non-degenerate. We fix a root vector $u_\alpha$ and denote by $p(\alpha)$ the parity of $u_\alpha$ for $\alpha \in \Delta$.

Let $V^k(\mathfrak{g})$ be the affine vertex superalgebra associated to $\mathfrak{g}$ at level $k \in \C$, which is generated by $u(z)$ ($u \in \mathfrak{g}$) whose parity is the same as $u$, satisfying that
\begin{align*}
[u_\lambda v] = [u, v] + k(u|v)\lambda,\quad
u, v \in \mathfrak{g}.
\end{align*}
Let $F(\mathfrak{g}_{\frac{1}{2}})$ be the neutral vertex superalgebra associated to $\mathfrak{g}_{\frac{1}{2}}$, which is strongly generated by $\phi_{\alpha}(z)$ ($\alpha \in \Delta_{\frac{1}{2}}$) whose parity is equal to $p(\alpha)$, satisfying that
\begin{align*}
[{\phi_\alpha}_\lambda \phi_\beta] = \chi(u_\alpha, u_\beta),\quad
\alpha, \beta \in \Delta_{\frac{1}{2}}.
\end{align*}
Let $F^\mathrm{ch}(\mathfrak{g}_{>0})$ be the charged fermion vertex superalgebra associated to $\mathfrak{g}_{>0}$, which is strongly generated by $\varphi_\alpha(z), \varphi^*_\alpha(z)$ ($\alpha \in \Delta_{>0}$) whose parities are equal to $p(\alpha) + \bar{1}$, satisfying that
\begin{align*}
[{\varphi_\alpha}_\lambda \varphi^*_\beta] = \delta_{\alpha, \beta},\quad
[{\varphi_\alpha}_\lambda \varphi_\beta] = [{\varphi^*_\alpha}_\lambda \varphi^*_\beta] = 0,\quad
\alpha, \beta \in \Delta_{>0}.
\end{align*}
Let $C^k(\mathfrak{g},f) = V^k(\mathfrak{g}) \otimes F(\mathfrak{g}_{\frac{1}{2}}) \otimes F^\mathrm{ch}(\mathfrak{g}_{>0})$ and $d$ be an odd element in $C^k(\mathfrak{g},f)$ defined by
\begin{align*}
d =& \sum_{\alpha \in \Delta_{>0}}(-1)^{p(\alpha)}u_\alpha\varphi^*_\alpha - \frac{1}{2}\sum_{\alpha, \beta, \gamma \in \Delta_{>0}}(-1)^{p(\alpha)p(\gamma)}c_{\alpha, \beta}^\gamma\NO{\varphi_\gamma\varphi^*_\alpha\varphi^*_\beta}\\
&+\sum_{\alpha \in \Delta_{\frac{1}{2}}}\phi_\alpha\varphi^*_\alpha + \sum_{\alpha \in \Delta_{>0}}\chi(u_\alpha)\varphi^*_\alpha.
\end{align*}
Then $(C^k(\mathfrak{g},f), d_{(0)})$ defines a cochain complex with respect to the charged degree: $\operatorname{charge}\varphi_\alpha = -\operatorname{charge}\varphi^*_\alpha = 1$ ($\alpha \in \Delta_{>0}$) and $\operatorname{charge}A=0$ for all $A \in V^k(\mathfrak{g})\otimes F(\mathfrak{g}_{\frac{1}{2}})$. The (affine) $\mathcal{W}$-algebra $\mathcal{W}^k(\mathfrak{g}, f)$ associated to $\mathfrak{g}$, $f$ at level $k$ is defined by
\begin{align*}
\mathcal{W}^k(\mathfrak{g}, f) = H(C^k(\mathfrak{g},f), d_{(0)}).
\end{align*}
Let $C^k(\mathfrak{g},f)_+$ be a subcomplex generated by $\phi_\alpha(z)$ ($\alpha \in \Delta_{\frac{1}{2}}$), $\varphi^*_\alpha(z)$ ($\alpha \in \Delta_{>0}$) and
\begin{align*}
J^u(z) = u(z) + \sum_{\alpha, \beta \in \Delta_{>0}}c_{\beta, u}^\alpha\NO{\varphi^*_\beta(z)\varphi_\alpha(z)},\quad
u \in \mathfrak{g}_{\leq 0}.
\end{align*}
Then we have \cite{KW04}
\begin{align*}
\mathcal{W}^k(\mathfrak{g}, f) = H(C^k(\mathfrak{g},f), d_{(0)}) = H^0(C^k(\mathfrak{g},f)_+, d_{(0)}).
\end{align*}
Thus, $\mathcal{W}^k(\mathfrak{g}, f)$ is a vertex subalgebra of $C^k(\mathfrak{g},f)_+$. Using the fact that
\begin{align*}
&[{J^u}_\lambda J^v] = J^{[u, v]} + \tau(u|v)\lambda,\quad
u, v \in \mathfrak{g}_{\leq0}\\
&\tau(u|v) = k(u|v) + \frac{1}{2}\kappa_{\mathfrak{g}}(u|v) - \frac{1}{2}\kappa_{\mathfrak{g}_0}(u|v),\quad
u, v \in \frak{g}_{\leq0},
\end{align*}
where $\kappa_{\mathfrak{g}}$ denotes the Killing form on $\mathfrak{g}$, it follows that the vertex algebra generated by $J^u(z)$ $(u \in \mathfrak{g}_{\leq0})$ is isomorphic to the affine vertex superalgebra associated to $\mathfrak{g}_{\leq0}$ and $\tau$, which we denote by $V^{\tau}(\mathfrak{g}_{\leq0})$. Therefore the homogeneous subspace of $C^k(\mathfrak{g},f)_+$ with charged degree $0$ is isomorphic to $V^{\tau}(\mathfrak{g}_{\leq0}) \otimes F(\mathfrak{g}_{\frac{1}{2}})$. The projection $\mathfrak{g}_{\leq 0} \twoheadrightarrow \mathfrak{g}_0$ induces a vertex superalgebra surjective homomorphism $V^{\tau}(\mathfrak{g}_{\leq0}) \otimes F(\mathfrak{g}_{\frac{1}{2}}) \twoheadrightarrow V^{\tau}(\mathfrak{g}_{0}) \otimes F(\mathfrak{g}_{\frac{1}{2}})$ so that we have
\begin{align*}
\Upsilon \colon \mathcal{W}^k(\mathfrak{g}, f) \rightarrow V^{\tau}(\mathfrak{g}_{0}) \otimes F(\mathfrak{g}_{\frac{1}{2}})
\end{align*}
by the restriction. The map $\Upsilon$ is called the Miura map and injective thanks to \cite{Frenkel, Arakawa17, Nakatsuka}.

\section{Finite $\mathcal{W}$-algebras}\label{sec:finite-W}
Recall the definitions of finite $\mathcal{W}$-algebras $U(\mathfrak{g}, f)$, following \cite{D3HK}. We introduce two definitions in \eqref{eq: finiteW-def1}, \eqref{eq: finiteW-def2} denoted by $U(\mathfrak{g}, f)_I$, $U(\mathfrak{g}, f)_{I\hspace{-.2em}I}$ respectively and prove the isomorphism $U(\mathfrak{g}, f)_I \simeq U(\mathfrak{g}, f)_{I\hspace{-.2em}I}$ in Theorem \ref{thm:D3HK}.
\smallskip

Let $\Phi$ be an associative $\C$-superalgebra generated by $\Phi_{\alpha}$ $(\alpha \in \Delta_{\frac{1}{2}})$ that has the same parity as $u_\alpha$, satisfying that
\begin{align*}
[\Phi_{\alpha}, \Phi_{\beta}] = \chi([u_\alpha, u_\beta]),\quad
\alpha, \beta \in \Delta_{\frac{1}{2}}.
\end{align*}
Here $[A,B]$ denotes $AB - (-1)^{p(A)\,p(B)}BA$. We extend the definition of $\Phi_\alpha$ for all $\alpha\in \Delta_{>0}$ by $\Phi_{\alpha} = 0$ for $\alpha \in \Delta_{\geq 1}$. Let $\Lambda(\mathfrak{g}_{>0})$ be the Clifford superalgebra associated to $\mathfrak{g}_{>0}$, which is an associative $\C$-superalgebra generated by $\psi_\alpha, \psi^*_\alpha$ $(\alpha \in \Delta_{>0})$ with the opposite parity to that of $u_\alpha$, satisfying that
\begin{align*}
[\psi_\alpha, \psi^*_\beta] = \delta_{\alpha, \beta},\quad
[\psi_\alpha, \psi_\beta] = [\psi^*_\alpha, \psi^*_\beta] = 0,\quad
\alpha, \beta \in \Delta_{>0}.  
\end{align*}
The Clifford superalgebra $\Lambda(\mathfrak{g}_{>0})$ has the charged degree defined by $\deg(\psi_\alpha) = 1 = -\deg(\psi^*_\alpha)$ for all $\alpha \in \Delta_{>0}$. Set
\begin{align*}
&C_I = U(\mathfrak{g})\otimes\Phi\otimes\Lambda(\mathfrak{g}_{>0}),\quad
d_I=\operatorname{ad}(Q),\\
&Q = \sum_{\alpha\in\Delta_{>0}}(-1)^{p(\alpha)}X_\alpha\psi_\alpha-\frac{1}{2}\sum_{\alpha,\beta,\gamma\in\Delta_{>0}}(-1)^{p(\alpha)p(\gamma)}c_{\alpha,\beta}^\gamma\psi_\gamma\psi^*_\alpha\ \psi^*_\beta,\\
&X_\alpha = u_\alpha + (-1)^{p(\alpha)}(\Phi_{\alpha} + \chi(u_\alpha)),\quad
\alpha \in \Delta_{>0},
\end{align*}
where $c_{\alpha, \beta}^\gamma$ is the structure constant defined by $[u_\alpha, u_\beta] = \sum_{\gamma \in \Delta_{>0}} c_{\alpha, \beta}^\gamma u_\gamma$. Then a pair $(C_I, d_I)$ forms a cochain complex with respect to the charged degree on $\Lambda(\mathfrak{g}_{>0})$ and the cohomology
\begin{align}\label{eq: finiteW-def1}
U(\mathfrak{g}, f)_{I} = H^\bullet(C_I, d_I)
\end{align}
has a structure of an associative $\C$-superalgebra inherited from that of $C_I$. Let
\begin{align*}
j^u = u + \sum_{\alpha, \beta \in \Delta_{>0}}c_{\beta, u}^\alpha\psi^*_\beta\ \psi_\alpha,\quad
u \in \mathfrak{g}.
\end{align*}
Then
\begin{align*}
\operatorname{ad}(Q) \cdot \psi_{\alpha} = j^{u_\alpha} +  (-1)^{p(\alpha)}(\Phi_{\alpha} + \chi(u_\alpha)) = X_\alpha + \sum_{\alpha, \beta\Delta_{>0}}c_{\beta, u}^\alpha\psi^*_\beta\ \psi_\alpha,\quad
\alpha \in \Delta_{>0}.
\end{align*}
Let $C_-$ be the subalgebra of $C_I$ generated by $\psi_\alpha$, $\operatorname{ad}(Q) \cdot \psi_{\alpha}$ $(\alpha \in \Delta_{>0})$ and $C_+$ be the subalgebra of $C_I$ generated by $j^u$ $(u \in \mathfrak{g}_{\leq0})$, $\Phi_\alpha$ $(\alpha \in \Delta_{\frac{1}{2}})$ and $\psi^*_\alpha$ $(\alpha \in \Delta_{>0})$. Then $(C_\pm,  d_I)$ form subcomplexes and $C_I \simeq C_- \otimes C_+$ as vector superspaces. Since $H(C_-, d_I)=\C$, we have
\begin{align*}
H(C_I, d_I) \simeq\ H(C_-, d_I) \otimes H(C_+, d_I) = H(C_+, d_I).
\end{align*}
Using the same argument as in \cite{KW04}, it follows that $H^n(C_+, d_I)=0$ for $n\neq0$. Therefore $U(\mathfrak{g}, f)_{I}$ is a subalgebra of $C^0_+$, which is generated by $j^u$ $(u \in \mathfrak{g}_{\leq0})$ and $\Phi_\alpha$ $(\alpha \in \Delta_{\frac{1}{2}})$. Since $[j^u, j^v] = j^{[u, v]}$ for $u, v \in \mathfrak{g}_{\leq0}$, there exists an isomorphism $C^0_+ \simeq U(\mathfrak{g}_{\leq0}) \otimes \Phi$ as associative $\C$-superalgebras. The projection $\mathfrak{g}_{\leq 0} \twoheadrightarrow \mathfrak{g}_0$ induces an associative $\C$-superalgebra surjective homomorphism $U(\mathfrak{g}_{\leq 0}) \otimes \Phi \twoheadrightarrow U(\mathfrak{g}_{0}) \otimes \Phi$ so that we have
\begin{align*}
\mu \colon U(\mathfrak{g}, f)_{I} \rightarrow U(\mathfrak{g}_{0}) \otimes \Phi
\end{align*}
by the restriction. The map $\mu$ is called the Miura map for the finite $\mathcal{W}$-algebras and injective by \cite{Lynch, Genra20, Nakatsuka}. Let $\C_{-\chi}$ be the one-dimensional $\mathfrak{g}_{\geq1}$-module defined by $\mathfrak{g}_{\geq1} \ni u \mapsto -\chi(u) \in \C$ and $M_{I\hspace{-.2em}I}$ be the induced left $\mathfrak{g}$-module
\begin{align*}
M_{I\hspace{-.2em}I} = \operatorname{Ind}^\mathfrak{g}_{\mathfrak{g}_{\geq1}}\C_{-\chi} = U(\mathfrak{g})\underset{U(\mathfrak{g}_{\geq1})}{\otimes}\C_{-\chi} \simeq U(\mathfrak{g})/I_{-\chi},
\end{align*}
where $I_{-\chi}$ is a left $U(\mathfrak{g})$-module generated by $u + \chi(u)$ for all $u \in \mathfrak{g}_{\geq1}$. Then $M_{I\hspace{-.2em}I}$ has a structure of the $\operatorname{ad}(\mathfrak{g}_{>0})$-module inherited from that of $U(\mathfrak{g})$. Set the $\operatorname{ad}(\mathfrak{g}_{>0})$-invariant subspace
\begin{align}\label{eq: finiteW-def2}
U(\mathfrak{g}, f)_{I\hspace{-.2em}I} = (M_{I\hspace{-.2em}I})^{\operatorname{ad}(\mathfrak{g}_{>0})}.
\end{align}
Then $U(\mathfrak{g}, f)_{I\hspace{-.2em}I}$ also has a structure of an associative $\C$-superalgebra inherited from that of $U(\mathfrak{g})$. We may also define $U(\mathfrak{g}, f)_{I\hspace{-.2em}I}$ as the Chevalley cohomology $H(\mathfrak{g}_{>0}, M_{I\hspace{-.2em}I})$ of the left $\mathfrak{g}_{>0}$-module $M_{I\hspace{-.2em}I}$:

\begin{lemma}[{\cite{GG, Nakatsuka}}]
\begin{align*}
H(\mathfrak{g}_{>0}, M_{I\hspace{-.2em}I}) = H^0(\mathfrak{g}_{>0}, M_{I\hspace{-.2em}I}) = (M_{I\hspace{-.2em}I})^{\operatorname{ad}(\mathfrak{g}_{>0})}.
\end{align*}
\begin{proof}
Though the assertion is proved in \cite{GG} for Lie algebras $\mathfrak{g}$, the same proof together with \cite[Corollary 2.6]{Nakatsuka} applies.
\end{proof}
\end{lemma}

\begin{theorem}[{\cite[Theorem A.6]{D3HK}}]\label{thm:D3HK}
There exists an isomorphism $U(\mathfrak{g}, f)_I \simeq U(\mathfrak{g}, f)_{I\hspace{-.2em}I}$ as associative $\C$-superalgebras.
\begin{proof}
Though the assertion is proved in \cite{D3HK} for Lie algebras $\mathfrak{g}$, the same proof applies as follows. Let $C_{I\hspace{-.2em}I} = \Lambda(\mathfrak{g}_{>0})_c \otimes M_{I\hspace{-.2em}I}$ be the Chevalley cohomology complex of the left $\mathfrak{g}_{>0}$-module $M_{I\hspace{-.2em}I}$, where $\Lambda(\mathfrak{g}_{>0})_c$ is the subalgebra of $\Lambda(\mathfrak{g}_{>0})$ generated by $\psi^*_\alpha$ for all $\alpha \in \Delta_{>0}$, and $d_{I\hspace{-.2em}I}$ be the derivation of the cochain complex $C_{I\hspace{-.2em}I}$. Let $U(\mathfrak{g}_{>0})_{-\chi} = U(\mathfrak{g}_{>0})\otimes\C_{-\chi}$ be a left $\mathfrak{g}_{\geq1}$-module defined by the diagonal action, where $U(\mathfrak{g}_{>0})$ is considered as a left $\mathfrak{g}_{\geq1}$-module by the left multiplication, and $M_{I\hspace{-.2em}I\hspace{-.2em}I}$ be the induced left $\mathfrak{g}$-module
\begin{align*}
M_{I\hspace{-.2em}I\hspace{-.2em}I}
= \operatorname{Ind}^\mathfrak{g}_{\mathfrak{g}_{\geq1}}U(\mathfrak{g}_{>0})_{-\chi}
= U(\mathfrak{g})\underset{U(\mathfrak{g}_{\geq1})}{\otimes}U(\mathfrak{g}_{>0})_{-\chi}.
\end{align*}
Let $\C_\chi$ be the one-dimensional $\mathfrak{g}_{\geq1}$-module defined by $\mathfrak{g}_{\geq1} \ni u \mapsto \chi(u) \in \C$ and $U(\mathfrak{g})_\chi = U(\mathfrak{g}) \otimes \C_\chi$ be a right $\mathfrak{g}_{\geq1}$-module defined by the diagonal action, where $U(\mathfrak{g})$ is considered as a right $\mathfrak{g}_{\geq1}$-module by the right multiplication. Then we have
\begin{align*}
M_{I\hspace{-.2em}I\hspace{-.2em}I} \simeq U(\mathfrak{g})_\chi\underset{U(\mathfrak{g}_{\geq1})}{\otimes}U(\mathfrak{g}_{>0})
\end{align*}
so that $M_{I\hspace{-.2em}I\hspace{-.2em}I}$ is a left $\mathfrak{g}$- right $\mathfrak{g}_{>0}$-bimodule. Note that there is an isomorphisms $\Lambda(\mathfrak{g}_{>0}) \simeq \Lambda(\mathfrak{g}_{>0})_h \otimes \Lambda(\mathfrak{g}_{>0})_c$ of vector superspaces, where $\Lambda(\mathfrak{g}_{>0})_h$ is the subalgebra of $\Lambda(\mathfrak{g}_{>0})$ generated by $\psi_\alpha$ for all $\alpha \in \Delta_{>0}$. Let $d_h$ be the derivation of the Chevalley homology complex $M_{I\hspace{-.2em}I\hspace{-.2em}I}\otimes\Lambda(\mathfrak{g}_{>0})_h$ of the right $\mathfrak{g}_{>0}$-module $M_{I\hspace{-.2em}I\hspace{-.2em}I}$. Then $M_{I\hspace{-.2em}I\hspace{-.2em}I}\otimes\Lambda(\mathfrak{g}_{>0})_h$ is clearly a left $\mathfrak{g}_{>0}$-module with respect to the adjoint $\mathfrak{g}_{>0}$-action. Now, let $\overline{d}_c$ be the derivation of the Chevalley cohomology complex $\Lambda(\mathfrak{g}_{>0})_c \otimes M_{I\hspace{-.2em}I\hspace{-.2em}I}\otimes\Lambda(\mathfrak{g}_{>0})_h$ of the left $\mathfrak{g}_{>0}$-module $M_{I\hspace{-.2em}I\hspace{-.2em}I}\otimes\Lambda(\mathfrak{g}_{>0})_h$. Then, as in \cite{D3HK}, we get a new cochain complex $(C_{I\hspace{-.2em}I\hspace{-.2em}I}, d_{I\hspace{-.2em}I\hspace{-.2em}I})$ defined by
\begin{align*}
C_{I\hspace{-.2em}I\hspace{-.2em}I} = \Lambda(\mathfrak{g}_{>0})_c \otimes M_{I\hspace{-.2em}I\hspace{-.2em}I}\otimes\Lambda(\mathfrak{g}_{>0})_h,\quad
d_{I\hspace{-.2em}I\hspace{-.2em}I} = d_c +(-1)^{\delta-1}\otimes d_h,
\end{align*}
where $\delta$ denotes the parity of the part of elements in $\Lambda(\mathfrak{g}_{>0})_c$. Then it is easy to check that the following linear map
\begin{multline*}
i_{I\hspace{-.2em}I\hspace{-.2em}I \rightarrow I} \colon C_{I\hspace{-.2em}I\hspace{-.2em}I} \ni\  \psi^*_{\beta_1}\cdots \psi^*_{\beta_i} \otimes (v_1 \cdots v_s \underset{U(\mathfrak{g}_{\geq1})}{\otimes}u_{\alpha_1} \cdots u_{\alpha_t}) \otimes \psi_{\gamma_1}\cdots \psi_{\gamma_j}\\
\mapsto\ \psi^*_{\beta_1}\cdots \psi^*_{\beta_i} \cdot v_1 \cdots v_s \cdot X_{\alpha_1} \cdots X_{\alpha_t} \cdot \psi_{\gamma_1}\cdots \psi_{\gamma_j} \in \overline{C}_I
\end{multline*}
with $v_1, \ldots, v_s \in \mathfrak{g}$, $\alpha_1,\ldots,\alpha_t, \beta_1, \ldots,\beta_i, \gamma_1, \ldots, \gamma_j \in \Delta_{>0}$ is well-defined and induces an isomorphism of complexes $(C_{I\hspace{-.2em}I\hspace{-.2em}I}, d_{I\hspace{-.2em}I\hspace{-.2em}I}) \rightarrow (C_I, d_I )$ since $i_{I\hspace{-.2em}I\hspace{-.2em}I \rightarrow I} \circ d_{I\hspace{-.2em}I\hspace{-.2em}I} = d_I \circ i_{I\hspace{-.2em}I\hspace{-.2em}I \rightarrow I}$. Now
\begin{align*}
H_n(C_{I\hspace{-.2em}I\hspace{-.2em}I}, d_h)
&= \Lambda(\mathfrak{g}_{>0})_c \otimes H_n\left(M_{I\hspace{-.2em}I\hspace{-.2em}I}\otimes\Lambda(\mathfrak{g}_{>0})_h, d_h\right)\\
&= \Lambda(\mathfrak{g}_{>0})_c \otimes U(\mathfrak{g})_\chi\underset{U(\mathfrak{g}_{\geq1})}{\otimes} H_n(\mathfrak{g}_{>0}, U(\mathfrak{g}_{>0}))\\
&= \delta_{n, 0}\ \Lambda(\mathfrak{g}_{>0})_c \otimes U(\mathfrak{g})_\chi\underset{U(\mathfrak{g}_{\geq1})}{\otimes}\C
\simeq \delta_{n,0}\ C_{I\hspace{-.2em}I}.
\end{align*}
Thus, since $d_c$ and $(-1)^{\delta-1}\otimes d_h$ commute, we have
\begin{align*}
H(C_{I\hspace{-.2em}I\hspace{-.2em}I}, d_{I\hspace{-.2em}I\hspace{-.2em}I})
\simeq H(H(C_{I\hspace{-.2em}I\hspace{-.2em}I},d_h), d_c)
\simeq H(C_{I\hspace{-.2em}I}, d_{I\hspace{-.2em}I}).
\end{align*}
The above argument together with the isomorphism $i_{I\hspace{-.2em}I\hspace{-.2em}I \rightarrow I}$ of complexes shows that $(C_I, d_I )$ and $(C_{I\hspace{-.2em}I}, d_{I\hspace{-.2em}I})$ are quasi-isomorphic via the following quasi-isomorphism
\begin{multline}\label{eq:D3HK-quasi}
i_{I \rightarrow I\hspace{-.2em}I}\colon C_I \ni\ \psi^*_{\beta_1}\cdots \psi^*_{\beta_i} \cdot v_1 \cdots v_s \cdot X_{\alpha_1} \cdots X_{\alpha_t} \cdot \psi_{\gamma_1}\cdots \psi_{\gamma_j}\\
\mapsto\ \delta_{t,0}\delta_{j,0}\ \psi^*_{\beta_1}\cdots \psi^*_{\beta_i} \cdot v_1 \cdots v_s \in C_{I\hspace{-.2em}I},
\end{multline}
which preserves the associative superalgebra structures on the cohomologies.
\end{proof}
\end{theorem}

\begin{definition}
The finite $\mathcal{W}$-algebra $U(\mathfrak{g}, f)$ associated to $\mathfrak{g}, f$ is defined to be the superalgebra $U(\mathfrak{g}, f)_{I}$, which is isomorphic to $U(\mathfrak{g}, f)_{I\hspace{-.2em}I}$ due to Theorem \ref{thm:D3HK}.
\end{definition}

\begin{remark}\label{rem:def-W}
The same results as Theorem \ref{thm:D3HK} for Poisson superalgebra versions has been studied in \cite{Suh16}. Also remark that our definitions of the finite $\mathcal{W}$-algebra $U(\mathfrak{g}, f)$ are not necessarily equivalent to the definitions in some literatures \cite{Poletaeva13, PS, ZS}. In fact, in case that $\mathfrak{g}=\mathfrak{osp}_{1|2n}$ and $f = f_\mathrm{prin}$ its principal nilpotent element, we have $\dim\mathfrak{g}_{\frac{1}{2}} = \dim\mathfrak{g}_{\frac{1}{2}, \bar{1}} = 1$ and thus $\mathfrak{g}_{\geq1} \subsetneq \mathfrak{g}_{>0}$. Then $U(\mathfrak{g}, f) \simeq U(\mathfrak{g}, f)_{I\hspace{-.2em}I} = (U(\mathfrak{g})/I_{-\chi})^{\operatorname{ad}(\mathfrak{g}_{>0})}$ is a proper subalgebra of $(U(\mathfrak{g})/I_{-\chi})^{\operatorname{ad}(\mathfrak{g}_{\geq 1})} = \End_{U(\mathfrak{g})}U(\mathfrak{g})/I_{-\chi}$.
\end{remark}


The vertex superalgebra $C^k(\mathfrak{g}, f)$ has a conformal vector $\omega$ if $k\neq-h^\vee$, which defines the conformal weights on $C^k(\mathfrak{g}, f)$ by $L_0$, where $\omega(z) = \sum_{n \in \Z}L_n z^{-n-2}$. See \cite{KRW} for the details. Then $H=L_0$ defines a Hamiltonian operator on $C^k(\mathfrak{g}, f)$, the vertex subalgebra $C^k(\mathfrak{g},f)_+$ and the corresponding $\mathcal{W}$-algebra $\mathcal{W}^k(\mathfrak{g}, f)$. Moreover the Hamiltonian operator $L_0$ is well-defined for all $k \in \C$. Recall that $\operatorname{Zhu}_H V$ is the $H$-twisted Zhu algebra of $V$, see Section \ref{sec:Gamma/Z}. Let $x \in \mathfrak{h}$ such that $[x ,u] = j u$ for $u \in \mathfrak{g}_j$. Then by \cite{Arakawa17,DK},
\begin{align}\label{eq:Zhu-explicit}
\operatorname{Zhu}_H C^k(\mathfrak{g},f)_+ \simeq C_+,\quad
J^u \mapsto j^u + \tau(x|u),\ 
\phi_\alpha \mapsto \Phi_\alpha,\ 
\varphi_\alpha^* \mapsto \psi^*_\alpha
\end{align}
for $u \in \mathfrak{g}_{\leq 0}$, $\alpha \in \Delta_{>0}$ and $\operatorname{Zhu}_H H^0(C^k(\mathfrak{g},f)_+,d_{(0)}) \simeq H^0(C_+, d_I)$ so that 
\begin{align}\label{eq:tw-Zhu}
\operatorname{Zhu}_H \mathcal{W}^k(\mathfrak{g}, f) \simeq U(\mathfrak{g}, f).
\end{align}
Let $V_1, V_2$ be any $\frac{1}{2}\Z_{\geq0}$-graded vertex superalgebras with the Hamiltonian operators and $g \colon V_1\rightarrow V_2$ any vertex superalgebra homomorphism preserving the conformal weights. Since $g(V_1\circ V_1)=g(V_1)\circ g(V_1)\subset V_2\circ V_2$, the map $g$ induces an algebra homomorphism
\begin{align*}
\operatorname{Zhu}_H (g)\colon\operatorname{Zhu}_H V_1\rightarrow\operatorname{Zhu}_H V_2.
\end{align*}
Apply for $g = \Upsilon$. Then we get
\begin{align*}
\operatorname{Zhu}_H(\Upsilon) = \mu
\end{align*}
by construction.

\section{Principal $\mathcal{W}$-algebras of $\mathfrak{osp}_{1|2n}$}\label{sec:prin-W}

Consider the case that
\begin{align*}
\mathfrak{g} = \mathfrak{osp}_{1|2n}
= \left\{ u =
\left(
\begin{array}{c|cc}
0 & {}^t y & -{}^t x \\
\hline %
x & a & b \\
y & c & -{}^t a
\end{array}
\right) \in \mathfrak{gl}_{1|2n}
\mathrel{} \middle| \mathrel{}
\begin{array}{l}
a, b, c \in \operatorname{Mat}_\C(n \times n),\\
x, y \in \operatorname{Mat}_\C(n \times 1),\\
b = {}^t b,\ c = {}^t c
\end{array}
\right\},
\end{align*}
where ${}^t A$ denotes the transpose of $A$. Let $\{e_{i, j}\}_{i, j \in I}$ be the standard basis of $\mathfrak{gl}_{1|2n}$ with the index set $I = \{ 0, 1, \ldots, n, -1, \ldots, -n\}$ and $h_i = e_{i, i} - e_{-i, -i}$ $(i = 1, \ldots, n)$. Then $\mathfrak{h} = \operatorname{Span}_\C\{h_i\}_{i=1}^n$ is a Cartan subalgebra of $\mathfrak{osp}_{1|2n}$. Define $\epsilon_i \in \mathfrak{h}^*$ by $\epsilon_i(h_j) = \delta_{i, j}$. Then $\Delta_+  = \{\epsilon_i, 2\epsilon_i\}_{i=1}^n \sqcup \{\epsilon_i-\epsilon_j,\epsilon_i+\epsilon_j\}_{1 \leq i<j \leq n}$ forms a set of positive roots with simple roots $\Pi = \{ \alpha_i\}_{i=1}$, $\alpha_i = \epsilon_i - \epsilon_{i+1}$ $(i=1,\ldots,n-1)$ and $\alpha_n = \epsilon_n$, and $\epsilon_1, \ldots, \epsilon_n$ are the (non-isotropic) odd roots in $\Delta_+$. Set $\Delta_- = -\Delta_+$ and $(u|v) = -\operatorname{str}(uv)$ for $u, v \in \mathfrak{osp}_{1|2n}$. We may identify $\mathfrak{h}^*$ with $\mathfrak{h}$ through $\nu \colon \mathfrak{h}^* \ni \lambda \mapsto \nu(\lambda) \in \mathfrak{h}$ defined by $\lambda(h) = (h|\nu(\lambda))$ for $h \in \mathfrak{h}$, which induces a non-degenerate bilinear form on $\mathfrak{h}^*$ by $(\lambda|\mu) = (\nu(\lambda)|\nu(\mu))$ so that $(\epsilon_i|\epsilon_j) = \delta_{i, j}/2$. Then $h_i$ corresponds to $2\epsilon_i = 2 \sum_{j=i}^n \alpha_j$ by $\nu$. We have
\begin{align*}
(\alpha_i|\alpha_i) = 1,\quad
(\alpha_i|\alpha_{i+1}) = -\frac{1}{2},\quad
i = 1, \ldots, n-1;\quad
(\alpha_n|\alpha_n) = \frac{1}{2}.
\end{align*}
Note that the dual Coxeter number of $\mathfrak{osp}_{1|2n}$ is equal to $n+\frac{1}{2}$. Let
\begin{align*}
f_\mathrm{prin} = \sum_{i=1}^{n-1} u_{-\alpha_i} + u_{-2\alpha_n}
\end{align*}
be a principal nilpotent element in the even part of $\mathfrak{osp}_{1|2n}$, where $u_\alpha$ denotes some root vector for $\alpha \in \Delta$. Then there exists a unique good grading on $\mathfrak{osp}_{1|2n}$ such that $\Pi_1 = \{\alpha_i\}_{i=1}^{n-1}$ and $\Pi_{\frac{1}{2}} = \{\alpha_n\}$. Thus
\begin{align*}
\mathfrak{g}_0 = \mathfrak{h},\ 
\mathfrak{g}_{>0} = \mathfrak{n} := \bigoplus_{\alpha \in \Delta_+}\mathfrak{g}_\alpha,\ 
\mathfrak{g}_{<0} = \mathfrak{n}_- := \bigoplus_{\alpha \in \Delta_-}\mathfrak{g}_\alpha.
\end{align*}
Let
\begin{align*}
\mathcal{W}^k(\mathfrak{osp}_{1|2n}) := \mathcal{W}^k(\mathfrak{osp}_{1|2n}, f_\mathrm{prin})
\end{align*}
be the principal $\mathcal{W}$-algebra of $\mathfrak{osp}_{1|2n}$ at level $k$. The Miura map for $\mathcal{W}^k(\mathfrak{osp}_{1|2n})$ is
\begin{align*}
\Upsilon\colon \mathcal{W}^k(\mathfrak{osp}_{1|2n}) \rightarrow \pi \otimes F,
\end{align*}
where $\pi$ is the Heisenberg vertex algebra generated by even fields $\alpha_i(z)$ $(i = 1, \ldots, n)$ satisfying that
\begin{align*}
[{\alpha_i}_\lambda \alpha_j] =  \left(k+n+\frac{1}{2}\right)(\alpha_i|\alpha_j)\lambda,\quad
i, j = 1,\ldots,n
\end{align*}
and $F$ is the free fermion vertex superalgebra generated by an odd field $\phi(z)$ satisfying that
\begin{align*}
[\phi_\lambda\phi] = 1.
\end{align*}
By \cite[Theorem 6.4]{Genra17}, $\mathcal{W}^k(\mathfrak{osp}_{1|2n})$ is strongly generated by $G, W_2, W_4, \ldots, W_{2n}$ for odd $G$ and even $W_2, W_4, \ldots, W_{2n}$ elements of conformal weights $n + \frac{1}{2}$ and $2, 4, \ldots, 2n$ such that
\begin{equation}\label{eq:G_W_formula}
\begin{array}{ll}
&\displaystyle \Upsilon(G)(z) =\ \NO{(2(k+n)\partial + h_1(z)) \cdots (2(k+n)\partial + h_n(z))\phi(z)},\\[2mm]
&\displaystyle \Upsilon(W_{2i})(z) \equiv \sum_{1\leq j_1<\cdots<j_{i}\leq n}\NO{h_{j_1}^2(z)\cdots h_{j_{i}}^2(z)}\quad\left(\operatorname{mod}\ C_2(\pi\otimes F)\right),\\[5mm]
&\displaystyle C_2(\pi\otimes F)=\{ A_{(-2)}B\mid A, B\in\pi\otimes F\}.
\end{array}
\end{equation}
and
\begin{align}\label{eq:G-W-relation}
[G_\lambda G] = W_{2n} + \sum_{i=1}^{n-1}\gamma_i\left(\frac{\lambda^{2i-1}}{(2i-1)!}W_{2n-2i+1}+\frac{\lambda^{2i}}{(2i)!}W_{2n-2i}\right) + \gamma_n \frac{\lambda^{2n}}{(2n)!}
\end{align}
for some $W_{2j+1} \in \mathcal{W}^k(\mathfrak{osp}_{1|2n})$, where
\begin{align*}
h_i(z) = 2\sum_{j=i}^n \alpha_j(z),\quad
\gamma_i = (-1)^i\prod_{j=1}^i\left( 2(2j-1)(k+n)-1 \right)\left( 4j(k+n)+1 \right),
\end{align*}
which satisfy that
\begin{align*}
[{h_i}_\lambda h_j] = (2k+2n+1)\delta_{i, j}\lambda,\quad
i, j=1,\ldots, n.
\end{align*}
If $k + n + \frac{1}{2} \neq 0$,
\begin{align*}
L = \frac{W_2}{2(2k+2n+1)}
\end{align*}
is a unique conformal vector of $\mathcal{W}^k(\mathfrak{osp}_{1|2n})$ with the central charge
\begin{align*}
c(k) = -\frac{(2n+1)(2(2n-1)(k+n)-1)(4n(k+n)+1)}{2(2k+2n+1)}.
\end{align*}

\section{Zhu algebras of $\mathcal{W}^k(\mathfrak{osp}_{1|2n})$}\label{sec:Zhu_prinW}

By \eqref{eq:tw-Zhu}, we have an isomorphism
\begin{align*}
\iota_1 \colon \operatorname{Zhu}_H \mathcal{W}^k(\mathfrak{osp}_{1|2n}) \xrightarrow{\simeq} U(\mathfrak{osp}_{1|2n}, f_\mathrm{prin}).
\end{align*}
Then $\iota_1$ is induced by \eqref{eq:Zhu-explicit}:
\begin{align*}
&\operatorname{Zhu}_H C^k(\mathfrak{osp}_{1|2n},f_\mathrm{prin}) \xrightarrow{\simeq} C_+,\\ 
&J^u \mapsto j^u + (2k+2n+1)(\rho_{\mathfrak{osp}}|u),\ 
\phi_\alpha \mapsto \Phi_\alpha,\ 
\varphi^*_\alpha \mapsto \psi^*_\alpha,
\end{align*}
where
\begin{align*}
\rho_{\mathfrak{osp}} = \frac{1}{2}\sum_{\alpha \in \Delta_+}(-1)^{p(\alpha)}\alpha.
\end{align*}
Let $\C[\mathfrak{h}^*] = U(\mathfrak{h})$ and set an isomorphism
\begin{align*}
&\iota_2 \colon \operatorname{Zhu}_H \pi\otimes \operatorname{Zhu}_H F \xrightarrow{\simeq} \C[\mathfrak{h}^*]\otimes\Phi,\\
&h_i \mapsto h_i + (2n-2i+1)\left(k+n+\frac{1}{2}\right),\quad
\phi_{\alpha_n} \mapsto \Phi_{\alpha_n}.
\end{align*}
Then we have a commutative diagram of Miura maps
\begin{equation*}
\SelectTips{cm}{}
\xymatrix@W15pt@H11pt@R12pt@C30pt{
\operatorname{Zhu}_H \mathcal{W}^k(\mathfrak{osp}_{1|2n})\ar[r]^-{\operatorname{Zhu}_H (\Upsilon)}\ar[d]_-{\iota_1}&\operatorname{Zhu}_H \pi\otimes \operatorname{Zhu}_H F\ar[d]^{\iota_2}\\
U(\mathfrak{osp}_{1|2n}, f_\mathrm{prin})\ar[r]^-{\mu}&\C[\mathfrak{h}^*]\otimes\Phi.
}
\end{equation*}

By \cite{DK}, $\operatorname{Zhu}_H\mathcal{W}^k(\mathfrak{osp}_{1|2n})$ has a PBW basis generated by $G, W_2, W_4, \ldots, W_{2n}$. By abuse of notation, we shall use the same notations for the generators of $U(\mathfrak{osp}_{1|2n}, f_\mathrm{prin})$ corresponding to $G, W_2, W_4, \ldots, W_{2n}$ by $\iota_1$.

\begin{lemma}\label{lem:mu(G)}
$\mu(G) = \left(h_1 + \rho_{\mathfrak{osp}}(h_1)\right)\left(h_2 + \rho_{\mathfrak{osp}}(h_2)\right) \cdots \left(h_n + \rho_{\mathfrak{osp}}(h_n)\right)\otimes\Phi_{\alpha_n}$.
\begin{proof}
We have
\begin{align*}
\Upsilon(G) =&\ \NO{(2(k+n)\partial + h_1) \cdots (2(k+n)\partial + h_n)\phi}\\
\equiv& (-(2n-1)(k+n) + h_1) * (-(2n-3)(k+n) + h_2)*\\
& \cdots * (-(k+n) + h_n)* \phi
\quad\left(\operatorname{mod}\ \mathcal{W}^k(\mathfrak{osp}_{1|2n}) \circ \mathcal{W}^k(\mathfrak{osp}_{1|2n}) \right).
\end{align*}
Thus
\begin{align*}
\mu(G)=& \iota_2\Bigl((-(2n-1)(k+n) + h_1) * (-(2n-3)(k+n) + h_2) *\\
&\quad\cdots * (-(k+n) + h_n)* \phi\Bigr)\\
=&\left(h_1 + n-1 + \frac{1}{2}\right)\left(h_2 + n-2+\frac{1}{2}\right) \cdots \left(h_n + \frac{1}{2}\right)\otimes\Phi_{\alpha_n}.
\end{align*}
Therefore the assertion follows from the fact that $\rho_{\mathfrak{osp}}(h_i) = n-i+\frac{1}{2}$.
\end{proof}
\end{lemma}

For a basic classical Lie superalgebra $\mathfrak{g}$ such that $\mathfrak{g}_{\bar{1}} \neq 0$, denote by
\begin{align*}
&Z(\mathfrak{g}) = \{ z \in U(\mathfrak{g}) \mid u z - (-1)^{p(u)p(z)}z u = 0\ \mathrm{for}\ \mathrm{all}\ u \in \mathfrak{g}\},\\
&\mathcal{A}(\mathfrak{g}) = \{ a \in U(\mathfrak{g}) \mid u a - (-1)^{p(u)(p(a)+\bar{1})}a u = 0\ \mathrm{for}\ \mathrm{all}\ u \in \mathfrak{g}\},\\
&\widetilde{Z}(\mathfrak{g}) = Z(\mathfrak{g}) \oplus \mathcal{A}(\mathfrak{g}),
\end{align*}
called the center, the anticenter and the ghost center of $U(\mathfrak{g})$, respectively due to \cite{Gorelik01}. Then the ghost center $\widetilde{Z}(\mathfrak{g})$ coincides with the center of $U(\mathfrak{g})_{\bar{0}}$ by \cite[Corollary 4.4.4]{Gorelik01}. In case that $\mathfrak{g}=\mathfrak{osp}_{1|2n}$, there exists $T \in U(\mathfrak{g})_{\bar{0}}$ \cite{ABF, Musson97, GL} such that
\begin{align*}
\mathcal{A}(\mathfrak{osp}_{1|2n}) = Z(\mathfrak{osp}_{1|2n}) T,\quad
(\sigma \circ \eta)(T) = h_1h_2\cdots h_n,
\end{align*}
where
\begin{align*}
\eta \colon U(\mathfrak{osp}_{1|2n}) \twoheadrightarrow U(\mathfrak{h}) = \C[\mathfrak{h}^*]
\end{align*}
is the projection induced by the decomposition $U(\mathfrak{osp}_{1|2n}) \simeq \mathfrak{n}_- U(\mathfrak{osp}_{1|2n}) \oplus U(\mathfrak{h}) \oplus U(\mathfrak{osp}_{1|2n}) \mathfrak{n}$ and $\sigma$ is an isomorphism defined by
\begin{align*}
\sigma \colon \C[\mathfrak{h}^*] \rightarrow \C[\mathfrak{h}^*],\quad
f \mapsto (\sigma(f) \colon \lambda \mapsto f(\lambda-\rho_{\mathfrak{osp}}))
\end{align*}
The element $T$ is called the Casimir ghost \cite{ABF} since $T^2 \in Z(\mathfrak{osp}_{1|2n})$ such that $(\sigma \circ \eta)(T^2) = h_1^2\cdots h_n^2$, and is studied for general $\mathfrak{g}$ in \cite{Gorelik01}. It is well-known \cite{Kac84, Gorelik04} that the restriction of $\sigma \circ \eta$ to $Z(\mathfrak{g})$ is injective and maps onto $\C[\mathfrak{h}^*]^W$, where $W$ is the Weyl group of $\mathfrak{sp}_{2n}$, called the Harish-Chandra homomorphism of $\mathfrak{osp}_{1|2n}$. Recall that
\begin{align*}
U(\mathfrak{osp}_{1|2n}, f_\mathrm{prin}) \simeq U(\mathfrak{osp}_{1|2n}, f_\mathrm{prin})_{I\hspace{-.2em}I} = (U(\mathfrak{osp}_{1|2n})/I_{-\chi})^{\operatorname{ad}\mathfrak{n}},
\end{align*}
where $I_{-\chi}$ is a left $U(\mathfrak{osp}_{1|2n})$-module generated by $u_\alpha + (f_\mathrm{prin}|u_\alpha)$ for all $\alpha \in \Delta_+\setminus\{\alpha_n\}$. Define the projections $q_1, q_2$ by
\begin{align*}
&q_1 \colon U(\mathfrak{osp}_{1|2n}) \twoheadrightarrow U(\mathfrak{osp}_{1|2n})/I_{-\chi},\\
&q_2 \colon U(\mathfrak{osp}_{1|2n})/I_{-\chi} \simeq \mathfrak{n}_-U(\mathfrak{osp}_{1|2n})/I_{-\chi} \oplus U(\mathfrak{h}) \oplus U(\mathfrak{h})u_{\alpha_n} \twoheadrightarrow U(\mathfrak{h}) \oplus U(\mathfrak{h})u_{\alpha_n}
\end{align*}
and a linear map $q_3$ by
\begin{align*}
q_3 \colon U(\mathfrak{h}) \oplus U(\mathfrak{h})u_{\alpha_n} \rightarrow \C[\mathfrak{h}^*] \otimes \Phi,\quad
(f_1,\mathrel{} f_2 \cdot u_{\alpha_n}) \mapsto f_1\otimes1+ f_2\otimes\Phi_{\alpha_n}.
\end{align*}
Then, using the quasi-isomorphism $i_{I \rightarrow I\hspace{-.2em}I}$ in \eqref{eq:D3HK-quasi}, the Miura map $\mu$ can be identified with the restriction of the composition map $q_3 \circ q_2$ to $U(\mathfrak{osp}_{1|2n}, f_\mathrm{prin})_{I\hspace{-.2em}I}$ since $u_{\alpha_n} = X_{\alpha_n} + \Phi_{\alpha_n}$.

\begin{lemma}\label{lem:p(Tu)}
$q_1(Tu_{\alpha_n})$ is the element of $U(\mathfrak{osp}_{1|2n}, f_\mathrm{prin})_{I\hspace{-.2em}I}$ corresponding to $G$.
\begin{proof}
First of all, we show that $q_1(Tu_{\alpha_n}) \in U(\mathfrak{osp}_{1|2n}, f_\mathrm{prin})_{I\hspace{-.2em}I}$. It is enough to show that $[u_\alpha, Tu_{\alpha_n}] \equiv 0$ $(\operatorname{mod}. I_{-\chi})$ for all $\alpha \in \Delta_+$. Let $\Delta_{+, \bar{i}} = \{ \alpha \in \Delta_+ \mid p(u_\alpha) = \bar{i}\}$. Since $[u_\alpha, T]=0$ for $\alpha \in \Delta_{+, \bar{0}}$, we have
\begin{align*}
[u_\alpha, Tu_{\alpha_n}] = T[u_\alpha, u_{\alpha_n}] \equiv 0\quad
(\operatorname{mod}. I_{-\chi}),\quad
\alpha \in \Delta_{+, \bar{0}}.
\end{align*}
Next, for $\alpha \in \Delta_{+, \bar{1}}\setminus\{\alpha_n\}$, since $u_\alpha T + Tu_\alpha = 0$, we also have
\begin{align*}
[u_\alpha, Tu_{\alpha_n}] = -T[u_\alpha, u_{\alpha_n}] + 2Tu_{\alpha_n}u_\alpha\equiv 0\quad
(\operatorname{mod}. I_{-\chi}),\quad
\alpha \in \Delta_{+, \bar{1}}\setminus\{\alpha_n\}.
\end{align*}
Finally, in case that $\alpha = \alpha_n$,
\begin{align*}
[u_{\alpha_n}, Tu_{\alpha_n}] = (u_{\alpha_n}T  + Tu_{\alpha_n})u_{\alpha_n} = 0.
\end{align*}
Therefore, $q_1(Tu_{\alpha_n})$ belongs to $U(\mathfrak{osp}_{1|2n}, f_\mathrm{prin})_{I\hspace{-.2em}I}$. Now $\mu = q_3 \circ q_2|_{U(\mathfrak{osp}_{1|2n}, f_\mathrm{prin})_{I\hspace{-.2em}I}}$ and by definition,
\begin{align*}
((\sigma \otimes 1) \circ \mu )(q_1(Tu_{\alpha_n})) &= ((\sigma \otimes 1) \circ q_3 \circ q_2 \circ q_1)(Tu_{\alpha_n})\\
&= (\sigma \circ \eta)(T) \otimes \Phi_{\alpha_n} = h_1\cdots h_n \otimes\Phi_{\alpha_n}.
\end{align*}
By Lemma \ref{lem:mu(G)}, $((\sigma \otimes 1) \circ \mu )(G) = h_1\cdots h_n \otimes\Phi_{\alpha_n}$. Since $(\sigma \otimes 1) \circ \mu$ is injective, we have $q_1(Tu_{\alpha_n}) = G$.
\end{proof}
\end{lemma}

\begin{theorem}\label{thm:Uosp=Z}
$U(\mathfrak{osp}_{1|2n}, f_\mathrm{prin})_{\bar{0}} \simeq Z(\mathfrak{osp}_{1|2n})$.
\begin{proof}
Since $U(\mathfrak{osp}_{1|2n}, f_\mathrm{prin})$ has a PBW basis generated by $G$, $W_2$, $W_4, \ldots, W_{2n}$ and $G$ is a unique odd generator, $U(\mathfrak{osp}_{1|2n}, f_\mathrm{prin})_{\bar{0}}$ has a PBW basis generated by $W_2$, $W_4, \ldots, W_{2n}$. Now $\Phi$ is superalgebra generated by $\Phi_{\alpha_n}$ with the relation $2 \Phi_{\alpha_n}^2 = \chi(u_{\alpha_n}, u_{\alpha_n})$. Thus $\mu$ maps $U(\mathfrak{osp}_{1|2n}, f_\mathrm{prin})_{\bar{0}}$ to $\C[\mathfrak{h}^*]$. By \eqref{eq:G_W_formula}, $\mu(W_{2i})$ for $i = 1, \ldots, n$ are algebraically independent in $\C[\mathfrak{h}^*]$ with degree $2i$ (but not necessary homogeneous). Now, by definition, $q_2 \circ q_1 =\eta$ on $Z(\mathfrak{osp}_{1|2n})$. Hence $q_2 \circ q_1 |_{Z(\mathfrak{osp}_{1|2n})}$ is injective. In particular, $q_1 |_{Z(\mathfrak{osp}_{1|2n})}$ is injective. Clearly, $q_1(Z(\mathfrak{osp}_{1|2n}))$ is $\operatorname{ad}\mathfrak{n}$-invariant. Thus, $U(\mathfrak{osp}_{1|2n}, f_\mathrm{prin})\simeq U(\mathfrak{osp}_{1|2n}, f_\mathrm{prin})_{I\hspace{-.2em}I} $ contains $Z(\mathfrak{osp}_{1|2n})$ through $q_1$. Moreover
\begin{align*}
\mu(Z(\mathfrak{osp}_{1|2n}))
= (q_3 \circ q_2 \circ q_1 )(Z(\mathfrak{osp}_{1|2n}))
= \eta(Z(\mathfrak{osp}_{1|2n}))
= \sigma^{-1}(\C[\mathfrak{h}^*]^W).
\end{align*}
Since $\C [ \mathfrak{h}^*]^W$ is a symmetric algebra of $h_1^2, \ldots, h_n^2$, $\mu(Z(\mathfrak{osp}_{1|2n}))$ must contain all $\mu(W_{2i})$ for $i = 1, \ldots, n$. Therefore
\begin{align*}
U(\mathfrak{osp}_{1|2n}, f_\mathrm{prin})_{\bar{0}} \simeq  Z(\mathfrak{osp}_{1|2n}).
\end{align*}
This completes the proof.
\end{proof}
\end{theorem}

\begin{corollary}\label{cor:evenZhu=Z}
$\left(\operatorname{Zhu}_H \mathcal{W}^k(\mathfrak{osp}_{1|2n})\right)_{\bar{0}} \simeq Z(\mathfrak{osp}_{1|2n})$.
\begin{proof}
The assertion is immediate from Theorem \ref{thm:Uosp=Z} and the fact that $\operatorname{Zhu}_H \mathcal{W}^k(\mathfrak{osp}_{1|2n})$ $\simeq$ $U(\mathfrak{osp}_{1|2n}, f_\mathrm{prin})$.
\end{proof}
\end{corollary}

Consider a linear isomorphism
\begin{align*}
\xi \colon \widetilde{Z}(\mathfrak{osp}_{1|2n}) = Z(\mathfrak{osp}_{1|2n}) \oplus \mathcal{A}(\mathfrak{osp}_{1|2n}) \xrightarrow{\simeq} Z(\mathfrak{osp}_{1|2n}) \oplus \mathcal{A}(\mathfrak{osp}_{1|2n})u_{\alpha_n}
\end{align*}
defined by $\xi(z, a) = (z, a\,u_{\alpha_n})$. Then by Lemma \ref{lem:p(Tu)} and the fact that $\mathcal{A}(\mathfrak{osp}_{1|2n}) = Z(\mathfrak{osp}_{1|2n}) T$, we have $(q_1 \circ \xi)(\widetilde{Z}(\mathfrak{osp}_{1|2n})) \subset U(\mathfrak{osp}_{1|2n}, f_\mathrm{prin})_{I\hspace{-.2em}I}$.

\begin{theorem}\label{thm:ghostZ_finiteW}
The map $q_1 \circ \xi \colon \widetilde{Z}(\mathfrak{osp}_{1|2n}) \rightarrow U(\mathfrak{osp}_{1|2n}, f_\mathrm{prin})$ is an isomorphisms of associative algebras.
\begin{proof}
By definition and Lemma \ref{lem:p(Tu)}, $(q_3 \circ q_2 \circ q_1 \circ \xi)(zT) = (q_3 \circ q_2 \circ q_1)(zTu_{\alpha_n}) = \eta(z)G$ for all $z \in Z(\mathfrak{osp}_{1|2n})$. Thus,
$q_3 \circ q_2 \circ q_1 \circ \xi|_{\mathcal{A}(\mathfrak{osp}_{1|2n})}$ is injective. In particular, $q_1 \circ \xi|_{\mathcal{A}(\mathfrak{osp}_{1|2n})}$ is injective. Using the fact that $U(\mathfrak{osp}_{1|2n}, f_\mathrm{prin})$ has a PBW basis generated by $G$, $W_2$, $W_4, \ldots, W_{2n}$ and Theorem \ref{thm:Uosp=Z}, it follows that $q_1 \circ \xi$ is a linear isomorphism. Now, we may suppose that $\chi(u_{\alpha_n}, u_{\alpha_n}) = 2$. Then $\Phi_{\alpha_n}^2 = 1$ so that $\mu(T^2) = \sigma^{-1}(h_1^2\cdots h_n^2) = \mu(G^2)$. Therefore $q_1 \circ \xi$ defines an isomorphisms of associative algebras. This completes the proof.
\end{proof}
\end{theorem}

Let $L(\lambda)$ be the simple highest weight $\mathfrak{osp}_{1|2n}$-module with the highest weight $\lambda$. Then there exists $\chi_\lambda \colon Z(\mathfrak{osp}_{1|2n}) \rightarrow \C$ such that $z$ acts on $\chi_\lambda(z)$ on $L(\lambda)$ for all $z \in Z(\mathfrak{osp}_{1|2N})$. The map $\chi_\lambda$ is called a central character of $\mathfrak{osp}_{1|2n}$ and is induced by $\eta$ and one-dimensional $\C[\mathfrak{h}^*]$-module $\C_\lambda$ defined by $f \mapsto f(\lambda)$. Using the Harish-Chandra homomorphism, it follows that $\chi_{\lambda_1} = \chi_{\lambda_2}$ if and only if $\lambda_2 = w(\lambda_1 + \rho_{\mathfrak{osp}}) - \rho_{\mathfrak{osp}}$ for some $w \in W$. Let
\begin{align*}
D = \{ \lambda \in \mathfrak{h}^* \mid \prod_{\alpha \in \Delta_{\bar{1}}}(\lambda+\rho_{\mathfrak{osp}}|\alpha) = 0\}.
\end{align*}
Denote by $\chi_\lambda \in D$ if $\lambda \in D$. Since $w(\Delta_{\bar{1}}) \subset \Delta_{\bar{1}}$ for all $w \in W$, we have $\lambda \in D \Rightarrow w(\lambda + \rho_{\mathfrak{osp}}) - \rho_{\mathfrak{osp}}\in D$ for any $w \in W$ so that $\chi_\lambda \in D$ is well-defined.

From now on, we will identify $\widetilde{Z}(\mathfrak{osp}_{1|2n})$ with $U(\mathfrak{osp}_{1|2n}, f_\mathrm{prin})$ by Theorem \ref{thm:ghostZ_finiteW}. Then $\widetilde{Z}(\mathfrak{osp}_{1|2n})$ is a superalgebra such that $\widetilde{Z}(\mathfrak{osp}_{1|2n})_{\bar{1}} = \mathcal{A}(\mathfrak{osp}_{1|2n})$. Let $E$ be a finite-dimensional $\Z_2$-graded simple $\widetilde{Z}(\mathfrak{osp}_{1|2n})$-module. Then $Z(\mathfrak{osp}_{1|2n})$ acts on $E$ as $\chi_\lambda$ for some $\lambda \in \mathfrak{h}^*$. For a non-zero parity-homogeneous element $v \in E$, $Tv$ has an opposite parity to $v$ such that $T^2 v = \chi_\lambda(T^2)v$. Recall that the set $\{h_1, \ldots, h_n\}$ is identified with $2\Delta_{+, \bar{1}}$ by $\mathfrak{h} \simeq \mathfrak{h}^*$. Then, using the fact that $\eta(T^2) = \sigma^{-1}(h_1^2\cdots h_n^2)$, it follows that
\begin{align*}
\chi_\lambda(T^2) = \prod_{i=1}^n\left((\lambda+\rho_{\mathfrak{osp}})(h_i)\right)^2 = \prod_{\alpha \in \Delta_{+, \bar{1}}}(\lambda + \rho_{\mathfrak{osp}}|2\alpha)^2.
\end{align*}
Hence $\chi_\lambda(T^2) = 0$ if and only if $\chi_\lambda \in D$. Since $E$ is simple, $E = \C v$ if $\chi_\lambda \in D$ and $E = \C v \oplus \C Tv$ if $\chi_\lambda \notin D$, which we denote by $E_{\chi_\lambda}$. Here we identify $E_{\chi_\lambda}$ with the parity change of $E_{\chi_\lambda}$ if $\chi_\lambda(T^2) = 0$. Therefore we obtain the following results:

\begin{proposition}\label{prop:rep_finiteW}
A finite-dimensional $\Z_2$-graded simple $U(\mathfrak{osp}_{1|2n}, f_\mathrm{prin})$-module is isomorphic to $E_{\chi_\lambda}$ for some $\lambda \in \mathfrak{h}^*$. In particular, there exists one-to-one correspondence between isomorphism classes (up to the parity change) of finite-dimensional $\Z_2$-graded simple $U(\mathfrak{osp}_{1|2n}, f_\mathrm{prin})$-modules and central characters of $\mathfrak{osp}_{1|2n}$.
\end{proposition}

\begin{corollary}\label{cor:Zhu_twisted}
There exists a bijective correspondence between central characters of $\mathfrak{osp}_{1|2n}$ and isomorphism classes (up to the parity change) of simple positive-energy Ramond-twisted $\mathcal{W}^k(\mathfrak{osp}_{1|2n})$-modules with finite-dimensional top spaces.
\begin{proof}
The assertion is immediate from $\operatorname{Zhu}_H \mathcal{W}^k(\mathfrak{osp}_{1|2n}) \simeq U(\mathfrak{osp}_{1|2n}, f_\mathrm{prin})$, Proposition \ref{prop:rep_finiteW} and \cite[Theorem 2.30]{DK}
\end{proof}
\end{corollary}

Corollary \ref{cor:Zhu_twisted} implies that dimensions of the top spaces $E_{\chi_\lambda}$ of simple positive-energy Ramond-twisted $\mathcal{W}^k(\mathfrak{osp}_{1|2n})$-modules are equal to $2$ if and only if $(\lambda+\rho_{\mathfrak{osp}}|\alpha) \neq 0$ for all $\alpha \in \Delta_{\bar{1}}$. We remark that this condition is equivalent to one that the annihilator of the Verma module $M(\lambda)$ is generated by its intersection with the center $Z(\mathfrak{osp}_{1|2n})$ by \cite{GL99}.

\bibliographystyle{halpha} 
\bibliography{refs} 

\end{document}